\documentclass[preprint,hidelinks,onefignum,onetabnum]{siamart220329}

\usepackage{amsmath}
\usepackage{amssymb}
\usepackage{amsfonts}
\usepackage{stmaryrd}
\usepackage{subcaption}

\usepackage{xspace}
\usepackage{xifthen}

\newcommand\mat[1]{\ensuremath\stackrel{\leftrightarrow}{\boldsymbol{#1}}}
\newcommand{\curl}[1][]{\ifthenelse{\isempty{#1}}{\operatorname{\nabla{\times}}}{\operatorname{\nabla_{#1}{\times}}}}
\newcommand{\diverge}[1][]{\ifthenelse{\isempty{#1}}{\operatorname{\nabla{\cdot}}}{\operatorname{\nabla_{#1}{\cdot}}}}

\usepackage{graphicx}

\usepackage{algorithm}
\usepackage{algpseudocode}
\usepackage{xspace}
\usepackage{xifthen}

\algnewcommand{\algorithmicgoto}{\textbf{go to}}%
\algnewcommand{\Goto}{\algorithmicgoto\xspace}%
\algnewcommand{\Label}{\State\unskip}

\makeatletter
\newenvironment{breakablealgorithm}
  {
   \begin{center}
     \refstepcounter{algorithm}
     \hrule height.8pt depth0pt \kern2pt
     \renewcommand{\caption}[2][\relax]{
       {\raggedright\textbf{\fname@algorithm~\thealgorithm} ##2\par}%
       \ifx\relax##1\relax 
         \addcontentsline{loa}{algorithm}{\protect\numberline{\thealgorithm}##2}%
       \else 
         \addcontentsline{loa}{algorithm}{\protect\numberline{\thealgorithm}##1}%
       \fi
       \kern2pt\hrule\kern2pt
     }
  }{
     \kern2pt\hrule\relax
   \end{center}
  }
\makeatother

\newsiamremark{remark}{Remark}
\newsiamremark{hypothesis}{Hypothesis}
\crefname{hypothesis}{Hypothesis}{Hypotheses}
\newsiamthm{claim}{Claim}

\headers{HDG methods for solving the two-fluid plasma model}{A. Ho and U. Shumlak}

\title{Hybridizable discontinuous Galerkin methods for solving the two-fluid plasma model\thanks{Submitted to the editors December 1, 2023.
\funding{The information, data, or work presented herein was funded in part by the Air Force Office of Scientific Research under award numbers FA9550-15-1-0271. This material is also based upon work supported by the National Science Foundation under Grant No. PHY-2108419.
}}}

\author{Andrew Ho\thanks{Computational Plasma Dynamics Lab, Aerospace and Energetics Research Program, University of Washington, Seattle, WA.}
\and Uri Shumlak\footnotemark[2]}

\usepackage{amsopn}

\DeclareMathOperator*{\sign}{sign}

\DeclareMathOperator{\grad}{\nabla}

\ifpdf
\hypersetup{
  pdftitle={Hybridizable discontinuous Galerkin methods for solving the two-fluid plasma model},
  pdfauthor={Andrew Ho and Uri Shumlak}
}
\fi

\begin{document}

\maketitle

\begin{abstract}
  The two-fluid plasma model has a wide range of timescales which must all be numerically resolved regardless of the timescale on which plasma dynamics occurs.
  The answer to solving numerically stiff systems is generally to utilize unconditionally stable implicit time advance methods.
  Hybridizable discontinuous Galerkin (HDG) methods have emerged as a powerful tool for solving stiff partial differential equations.
  The HDG framework combines the advantages of the discontinuous Galerkin (DG) method, such as high-order accuracy and flexibility in handling mixed hyperbolic/parabolic PDEs with the advantage of classical continuous finite element methods for constructing small numerically stable global systems which can be solved implicitly.
  In this research we quantify the numerical stability conditions for the two-fluid equations and demonstrate how HDG can be used to avoid the strict stability requirements while maintaining high order accurate results.
\end{abstract}



\section{Introduction}

The two-fluid plasma model is a very numerically stiff PDE system.
The model couples together the time evolution of an ion fluid, an electron fluid, and Maxwell's equations.
As a result the characteristic dynamics leads to timescales which can span two to four orders of magnitude, include both the speed of light and various electron wave modes along with ion fluid dynamics.
The large span of timescales makes the PDE system a prime candidate for being solved using implicit time stepping.
However, there are several challenges in actualizing an implicit numerical solution.
The system is a mixed advection-diffusion-reaction system consisting of six scalar fields and four vector fields with a nominal 18 components to be evolved in 3D.
The hyperbolic characteristics makes using standard continuous Galerkin (CG) methods challenging, and traditional discontinuous Galerkin (DG) methods are not typically amenable with implicit time stepping schemes\cite{peraire2008}.
One solution is the use of hybridizable discontinuous Galerkin (HDG) methods\cite{cockburn2009}.
These methods are designed to allow the implicit DG discretization to be partitioned into three distinct phases, reducing the largest global system solve size as well as producing a scheme which maintains local compactness.
Nguyen and Peraire\cite{nguyen2012} have demonstrated the applicability of HDG to a wide variety of linear and non-linear PDEs, including Stoke's equations, incompressible Navier-Stokes equations, and the compressible Euler equations. They demonstrated that the method can be formulated to achieve optimal $k+1$ orders of convergence, where $k$ is the polynomial basis order.
Lee et al.\cite{lee2018} demonstrated that the linearized incompressible resistive MHD equations can be solved using HDG, opening the door to modeling basic plasma systems using HDG.
In this paper we present techniques and tools for expanding the applicability of HDG to large coupled PDE systems such as the two-fluid plasma model using HDG.
\Cref{sec:tf_model} introduces the two-fluid plasma system to be solved, and discusses challenges associated with solving this particular system using explicit Runge-Kutta discontinuous Galerkin (RKDG).
\Cref{sec:hdg} presents the HDG discretization method and practical implementation details required for solving complex PDE systems.
\Cref{sec:convergence} presents the verification of the implementation by checking the convergence rates on linear systems.
Finally \cref{sec:plasma_tests} presents results of solving the two-fluid plasma using HDG.

\section{The two-fluid plasma model}
\label{sec:tf_model}

The two-fluid plasma model is a subset of the general multi-fluid 5N-moment plasma model\cite{shumlak2011}. This model treats each constituent particle species type as their own fluid species. Coupling between the species occurs via local inter-species collisions and long-range electromagnetic field coupling described using Maxwell's Equations.
For the two-fluid model a single positively charge ion species is coupled to an electron species.
The exact form being considered for this paper is given by
\begin{gather}
  \partial_{t} \rho_{\alpha} + \diverge \vec{p}_{\alpha} = 0\\
  \partial_{t} \vec{p}_{\alpha} + \diverge 
  \left(
    \vec{p}_{\alpha} \vec{u}_{\alpha} + P_{\alpha} \mat{I} -  \mu_{c} \sqrt{m_{\alpha}} \mat{W}_{\alpha}
  \right)
  = R_{c} \rho_{\alpha} (\vec{u}_{\beta} - \vec{u}_{\alpha}) + \frac{Z_{\alpha} L}{m_{\alpha} \delta_{p}} (\rho_{\alpha} \vec{E} + \vec{p}_{\alpha} \times \vec{B})\\
  \partial_{t} U_{\alpha}
  + \diverge 
  \left(
    U_{\alpha} \vec{u}_{\alpha} - \frac{\kappa_{c}}{\sqrt{m_{\alpha}}} \vec{h}_{\alpha}
  \right)
  + (P_{\alpha} \mat{I} - \mu_{c} \sqrt{m_{\alpha}} \mat{W}_{\alpha}) : \grad \vec{u}_{\alpha} = Q_{c} (T_{\beta} - T_{\alpha})\\
  \mat{W}_{\alpha} = \grad \vec{u}_{\alpha} + \grad^{T} \vec{u}_{\alpha} - \frac{2}{3} \mat{I} \diverge \vec{u}_{\alpha}\\
  \vec{h}_{\alpha} = \grad T_{\alpha}
\end{gather}
where $\rho_{\alpha}$, $\vec{p}_{\alpha}$, $U_{\alpha}$ are the mass, momentum, and internal energy densities of species $\alpha$ and $\mu_{c}$, $\kappa_{c}$, $R_{c}$, $Q_{c}$ are constants for the viscosity, thermal diffusivity, interspecies momentum transfer, and interspecies heat transfer coefficients respectively.
$m_{\alpha}$, $Z_{\alpha}$, and $\gamma$ describe the mass, charge state, and heat capacity ratio of species $\alpha$.
For the two-fluid model $\alpha$ is either $i$ for ions or $e$ for electrons, while $\beta$ is the other species.
The following additional substitutions are used to simplify the notation:
\begin{gather}
  \vec{u}_{\alpha} = \frac{\vec{p}_{\alpha}}{\rho_{\alpha}}\\
  P_{\alpha} = U_{\alpha} (\gamma - 1)\\
  T_{\alpha} = \frac{m_{\alpha} P_{\alpha}}{\rho_{\alpha}}\\
  \rho_{c} = \frac{Z_{i}}{m_{i}} \rho_{i} + \frac{Z_{e}}{m_{e}} \rho_{e}\\
  \vec{j} = \frac{Z_{i}}{m_{i}} \vec{p}_{i} + \frac{Z_{e}}{m_{e}} \vec{p}_{e}
\end{gather}
The plasma fluid equations are then coupled to EM fields via Maxwell's equations.
These are evolved using 
\begin{gather}
  \partial_{t} \vec{E} + \frac{c_{0}^{2} L}{V_{A}^{2} \delta_{p}} \vec{j}
  - \frac{c_{0}^{2}}{V_{A}^{2}} \curl \vec{B} - \grad \theta = 0\\
  \partial_{t} \vec{B} + \curl \vec{E} - \grad \psi = 0\\
  \frac{1}{c_{h}^{2}} \partial_{t} \psi + c_{p} \psi = \diverge \vec{B}\\
  \frac{1}{c_{h}^{2}} \partial_{t} \theta + c_{p} \theta = 
  \left(
    \diverge \vec{E} - \frac{c_{0}^{2} L}{V_{A}^{2} \delta_{p}} \rho_{c}
  \right)
\end{gather}
where $\vec{E}$ and $\vec{B}$ are the electric and magnetic fields respectively, $c_{0}/V_{A}$ is the speed of light normalized by the nominal Alfven speed $V_{A} = \sqrt{B_{0}^{2} / (\mu_{0} m_{p} n_{0})}$, and $\delta_{p}/L$ is the plasma skin depth $\delta_{p} = m_{p} V_{A} / (q_{e} B_{0})$ normalized by the characteristic length $L$.
This formulation includes a mixed hyperbolic-parabolic divergence cleaning operator to handle divergence errors which may be generated by the two-fluid plasma coupling\cite{munz2000,dedner2002}.

In total this system solves 42 coupled PDE equations and incorporates a simplified model of viscosity, heat flux, inter-species collision effects, and mixed hyperbolic-parabolic divergence cleaning.
In addition to the stiffness introduced by coupling parabolic diffusion and hyperbolic advection terms, the two-fluid plasma equations has characteristics which include\cite{shumlak2003} the speed of light in a vacuum, upper hybrid plasma oscillation frequency, along with the ion sound speed.

To quantify the stiffness of solving the two-fluid equation with explicit Runge-Kutta discontinuous Galerkin, a Von Neumann linear stability analysis\cite{charney1950} is performed on the ideal two-fluid equations. This simplified equation set has similar stability properties to the above collisional two-fluid equation set in inviscid regimes.
\begin{gather}
  \partial_{t} \rho_{\alpha} + \diverge \vec{p}_{\alpha} = 0\\
  \partial_{t} \vec{p}_{\alpha} + \diverge 
  \left(
    \vec{p}_{\alpha} \vec{u}_{\alpha} + P_{\alpha} \mat{I}
  \right)
  = \frac{Z_{\alpha} L}{m_{\alpha} \delta_{p}} (\rho_{\alpha} \vec{E} + \vec{p}_{\alpha} \times \vec{B})\\
  \partial_{t} e_{\alpha} + \diverge
  \left(
    (e_{\alpha} + P_{\alpha}) \vec{u}_{\alpha}
  \right) = \frac{Z_{\alpha} L}{m_{\alpha} \delta_{p}} \vec{p}_{\alpha} \cdot \vec{E}\\
  \partial_{t} \vec{E} + \frac{c_{0}^{2} L}{V_{A}^{2} \delta_{p}} \vec{j} - \frac{c_{0}^{2}}{V_{A}^{2}} \curl \vec{B} = 0\\
  \partial_{t} \vec{B} + \curl \vec{E} = 0
\end{gather}
where $e_{\alpha}$ is the total energy density of the fluid, related to the fluid internal energy density by
\begin{gather}
  e_{\alpha} = U_{\alpha} + \frac{1}{2} \vec{p}_{\alpha} \cdot \vec{u}_{\alpha}
\end{gather}
The ideal two-fluid plasma model can be written as a non-linear advection-reaction PDE of the form
\begin{gather}
  \partial_{t} \vec{q} + \diverge \mat{F}(\vec{q}) = \vec{S}(\vec{q})
\end{gather}
Rewriting using Jacobians and considering one spatial dimension,
\begin{gather}
  \partial_{t} \vec{q} + \mat{J}_{F} \cdot \partial_{x} \vec{q} = \mat{J}_{S} \cdot \vec{q}\\
  \mat{J}_{A} = \partial_{\vec{q}} \mat{F}(\vec{q})\\
  \mat{J}_{S} = \partial_{\vec{q}} \vec{S}(\vec{q})
\end{gather}
For the ideal two-fluid equations these can be expressed as block matrices such that
\begin{gather}
  \vec{q}^{T} =
  \begin{bmatrix}
    \rho_{i} & \vec{p}_{i} & e_{i} & \rho_{e} & \vec{p}_{e} & e_{e} & \vec{E} & \vec{B}
  \end{bmatrix}\\
  \mat{J}_{A} =
  \begin{bmatrix}
    \mat{J}_{A,i} & \mat{0} & \mat{0}\\
    \mat{0} & \mat{J}_{A,e} & \mat{0}\\
    \mat{0} & \mat{0} & \mat{J}_{A,EM}
  \end{bmatrix}\\
  \mat{J}_{S} =
  \begin{bmatrix}
    \mat{J}_{S,i,i} & \mat{0} & \mat{J}_{S,i,EM}\\
    \mat{0} & \mat{J}_{S,e,e} & \mat{J}_{S,e,EM}\\
    \mat{J}_{S,EM,i} & \mat{J}_{S,EM,e} & \mat{0}
  \end{bmatrix}
\end{gather}
This is then linearized by treating $\mat{J}_{A}$ and $\mat{J}_{S}$ as locally constant-coefficient matrices.
Then for each degree of freedom $q_{a} \in \vec{q}$ in an element define a perturbation
\begin{gather}
  q_{a}(x, t) = G_{a}(t) e^{j k x}
\end{gather}
where $j$ is the imaginary unit, $k$ is the perturbation wave number, and $G_{a} \in \vec{G}$ contains the temporal evolution dependence of $\vec{q}$.
Applying this to the two-fluid DG weak form with a Rusanov flux, we arrive at a semi-discrete operator $\mat{D}$ such that
\begin{gather}
  \partial_{t} \vec{G} = \mat{D} \vec{G}\\
  \begin{aligned}
    \mat{D}_{a + b N, c + d N}
    =&
       J_{A,b,d} A_{a,b} + \delta_{a,c} J_{S,b,d}
  + \delta_{b,d} \frac{\tau_{b}}{2} M^{-1}_{a,g} (L^{-}_{g} (L^{+}_{c} e^{-j k h} - L^{-}_{c}) - L^{+}_{g} (L^{+}_{c} - L^{-}_{c} e^{j k h}))\\
  &
  + \frac{J_{A,b,d}}{2} M^{-1}_{a,g} ( L^{-}_{g} (L^{-}_{c} + L^{+}_{c} e^{-j k h})
  - L^{+}_{g} (L^{+}_{c} + L^{-}_{c} e^{j k h})
  )
  \end{aligned}
\end{gather}
where
\begin{gather}
  M_{a,b} = \frac{h}{2} \int^{1}_{-1} \phi_{a}(\xi) \phi_{b}(\xi) d \xi\\
  A_{a,b} = M^{-1}_{a,c} \int^{1}_{-1} \phi_{b}(\xi) \partial_{\xi} \phi_{c}(\xi) d\xi\\
  L^{-}_{a} = \phi_{a}(-1)\\
  L^{+}_{a} = \phi_{a}(1)
\end{gather}
are the mass matrix, advection matrix, and lift operator respectively\cite{hesthaven2008}, $h$ is the element length, and $\vec{\tau}$ are the Rusanov fastest wave speeds.
For the two-fluid equations $\mat{J}_{A}$ is block diagonal and thus nominally separated into three independent hyperbolic systems, corresponding to the Euler equations for ions and electrons plus Maxwell's equations.
As such we will pick one Rusanov wave speed $\tau_{b} \in \vec{\tau}$ per block for a total of three unique values given by
\begin{gather}
  \tau_{i} = 
  \left|u_{x,i}\right| + \sqrt{\frac{\gamma P_{i}}{\rho_{i}}}
  \\
  \tau_{e} = 
  \left|u_{x,e}\right| + \sqrt{\frac{\gamma P_{e}}{\rho_{e}}}
  \\
  \tau_{EM} = \frac{c_{0}}{V_{A}}
\end{gather}
A numerical scheme is considered stable if all the eigenvalues of $\mat{D}\Delta t$ are within the region of absolute stability of the temporal discretization, for all $k h \in [0, 2 \pi]$.

\begin{figure}[H]
  \centering{}
  \includegraphics[width=.5\textwidth]{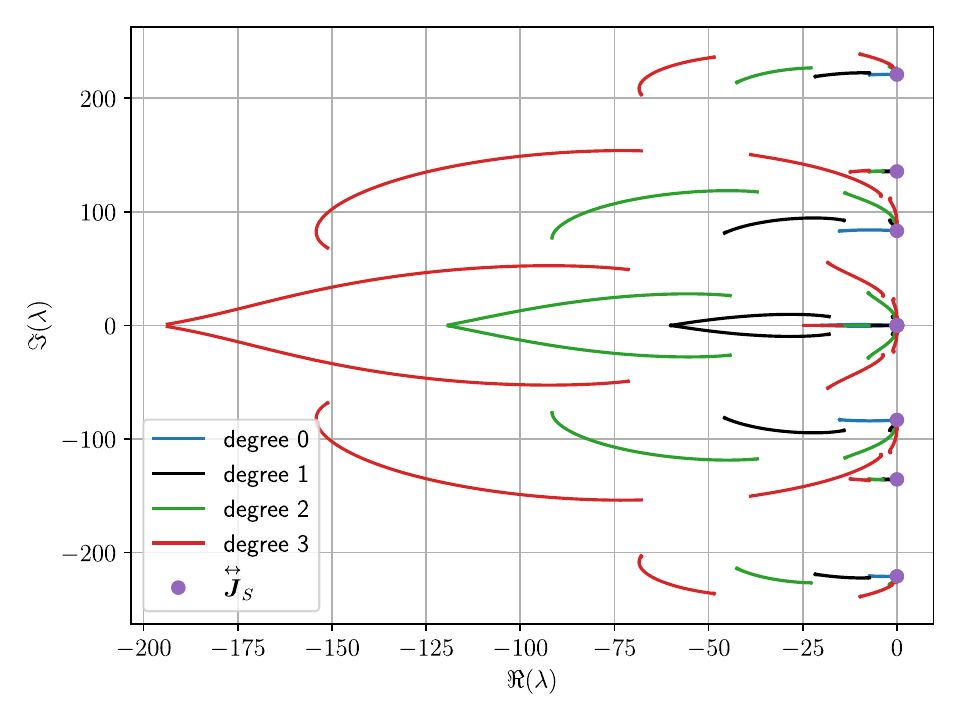}
  \caption{Example RKDG two-fluid linear stability eigencurves for $h = 10^{-1}$. Eigenvalues of $\mat{J}_{S}$ are plotted as points, and are guaranteed to be pure imaginary. The DG advection stability requirements of $\mat{J}_{A}$ grows quadratically with polynomial order while the eigenvalues of $\mat{J}_{S}$ are constant with respect to polynomial order leading to a transitions from $\mat{J}_{S}$ stability limited to $\mat{J}_{A}$ stability limited by increasing the polynomial order from zero to three\cite{warburton2008,chalmers2020}.}
  \label{fig:tf_eigs_dx_10}
\end{figure}

Consider the conditions where $\vec{u} = \vec{j} = 0$ and
\begin{gather}
  \begin{bmatrix}
    n_{i}\\
    U_{i}\\
    n_{e}\\
    U_{e}\\
    B_{x}\\
    B_{y}
  \end{bmatrix} =
                  \begin{bmatrix}
                    10\\
                    7.5 \times 10^{-5}\\
                    10\\
                    7.5 \times 10^{-5}\\
                    7.5 \times 10^{-2}\\
                    5 \times 10^{-5}
                  \end{bmatrix}
                  \\
                \begin{aligned}
                  \frac{\delta_{p}}{L} &= 1 & \gamma &= \frac{5}{3} & m_{i} &= 1 & m_{e} &= \frac{1}{1836} & \frac{c_{0}}{V_{A}} &= 1 & Z_{i} &= -Z_{e} = 1
  \end{aligned}
\end{gather}
\Cref{fig:tf_eigs_dx_10} shows the eigencurves which specify the eigenvalues of $\mat{D}$ as a parametric function of $kh$ for the above specified conditions.
There are a few notable features of the stability analysis to highlight:
\begin{enumerate}
\item There are ``pure'' hyperbolic advection modes, which originate from $\lambda = 0$ and have a negative real extent dictated by the DG discretization of $\mat{J}_{A}$.
  As a result the speed of light, which is the wave speed of Maxwell's equations, must be resolved within the temporal discretizations region of absolute stability.
\item The DG discretization eigencurves are known to grow quadratically with polynomial order\cite{warburton2008,chalmers2020}, and that remains true for the ideal two-fluid equations.
\item The eigencurves intercept the pure imaginary axis at the eigenvalues of $\mat{J}_{S}$. These are highlighted by the purple data points in \cref{fig:tf_eigs_dx_10}. As a result, temporal discretizations which are not stable for $\Im(\lambda) \neq 0$ such as forward Euler are unconditionally unstable.
\end{enumerate}
Having to resolve the speed of light is problematic when the timescales of interest are on the bulk plasma motion scale, which is on the order of the ion sound speed.
Similarly, it is possible that resolving the purely oscillatory modes from $\mat{J}_{S}$ will be the dominant factor on numerical stability.
As such, we examine these eigenvalues in greater detail.
The eigenvalues of $\mat{J}_{S}$ are given by
\begin{gather}
  \vec{\lambda}_{S} \in
  \begin{bmatrix}
    0\\
    \pm (\omega_{p} \tau) \sqrt{\omega_{p,i}^{2} + \omega_{p,e}^{2}}\\
    \pm 
    \left(
    \frac{L}{\delta_{p}}
    \right) \sqrt{z}
  \end{bmatrix}
\end{gather}
where
\begin{gather}
  \omega_{p} \tau = \frac{c_{0}}{V_{A}} \frac{L}{\delta_{p}}\\
  \omega_{p,\alpha} = \sqrt{\frac{n_{\alpha} Z_{\alpha}^{2}}{m_{\alpha}}}\\
  \omega_{c,\alpha} = \sqrt{\frac{Z_{\alpha}^{2}}{m_{\alpha}^{2}} \vec{B} \cdot \vec{B}}
\end{gather}
and $z$ are the roots of the polynomial
\begin{gather}
  \begin{aligned}
    0 =& z^{3}
    + \left(
    \omega_{c,i}^{2} + \omega_{c,e}^{2} + 2 (\omega_{p,i}^{2} + \omega_{p,e}^{2}) \frac{c_{0}^{2}}{V_{A}^{2}}
    \right) z^{2}\\
    &+
    \left(
    \omega_{c,i}^{2} \omega_{c,e}^{2} + 2
    \left(
    \omega_{p,i}^{2} \omega_{c,e}^{2} + \omega_{p,e}^{2} \omega_{c,i}^{2}
    \right) \frac{c_{0}^{2}}{V_{A}^{2}}
    + 
    \left(
    \omega_{p,i}^{2} + \omega_{p,e}^{2}
    \right)^{2} \frac{c_{0}^{4}}{V_{A}^{4}}
    \right) z\\
    &+ 
    \left(
    \omega_{p,e}^{2} \omega_{c,i} + \sign(Z_{i} Z_{e}) \omega_{p,i}^{2} \omega_{c,e}
    \right) \frac{c_{0}^{4}}{V_{A}^{4}}
  \end{aligned}
\end{gather}
We only consider cases where a quasi-neutral plasma could be formed, where is $Z_{i} Z_{e} < 0$.
There is no requirement that $Z_{i} > 0$ and $Z_{e} < 0$ or that $Z_{i} = - Z_{e}$, so long as the the two species are oppositely charged.
This allows the analysis to apply to positron/antiproton plasmas as well as multiply charged ion/electron plasmas.
We can prove the roots $z$ are non-positive real in two ways, via physics arguments or rigorous mathematics.
The physics argument which requires the roots to be non-positive real is that any eigenvalue $\lambda_{S}$ which is not pure imaginary (or zero) would result in a lack of conservation as the sources would have either an exponentially growing or decaying mode.
To prove this mathematically, consider each requirement separately: first prove the roots must be real, and separately prove that the roots must be in the negative half complex plane.
The first can be proven by showing the discriminant is non-negative.
Define
\begin{gather}
  \Gamma_{a} = \frac{\omega_{p,e}^{2}}{\omega_{c,i} (\omega_{c,i} + \omega_{c,e})} \frac{c_{0}^{2}}{V_{A}^{2}}\\
  \Gamma_{b} = \frac{\omega_{p,i}^{2}}{\omega_{c,e} (\omega_{c,i} + \omega_{c,e})} \frac{c_{0}^{2}}{V_{A}^{2}}
\end{gather}
Then re-arrange the discriminant $\Delta$ to be
\begin{gather}
  \begin{aligned}
    \Delta =&
             \left(
             \frac{c_{0}^{2}}{V_{A}^{2}} 
             \left(
             \omega_{p,i}^{2} \omega_{c,i} - \omega_{p,e}^{2} \omega_{c,e}
             \right)
             + \omega_{c,i} \omega_{c,e} (\omega_{c,i} - \omega_{c,e})
             \right)^{2}\\
             &\left(
             \begin{aligned}
               & 4 (\omega_{c,i} + \omega_{c,e})^{3}
                 \left(
                 \Gamma_{a} \omega_{c,i}^{3} (\Gamma_{a} - 1)^{2} + \Gamma_{b} \omega_{c,e}^{3} (\Gamma_{b} - 1)^{2}
                 \right)\\
               & + \omega_{c,i}^{2} \omega_{c,e}^{2} (\omega_{c,i} + \omega_{c,e})^{2}
                 \left(
                 (\Gamma_{b} - 1)^{2} + (\Gamma_{a} - 1)^{2}
                 \right)\\
               & + \Gamma_{a} \Gamma_{b} \omega_{c,i} \omega_{c,e} (19 (\omega_{c,i} + \omega_{c,e})^{4} + 12 (\omega_{c,i} + \omega_{c,e})^{3} (\Gamma_{a} \omega_{c,i} + \Gamma_{b} \omega_{c,e}) + (\omega_{c,i}^{2} - \omega_{c,e}^{2})^{2})\\
               & - \omega_{c,i}^{2} \omega_{c,e}^{2} (\omega_{c,i} + \omega_{c,e})^{2}
             \end{aligned}
             \right)
  \end{aligned}
\end{gather}
$\Delta$ is strictly non-negative if
\begin{gather}
  \Gamma_{a} \Gamma_{b} (19 + 12 (\Gamma_{a} + \Gamma_{b})) + (\Gamma_{a} - 1)^{2} + (\Gamma_{b} - 1) \ge 1
\end{gather}
At the equality point
\begin{gather}
  \Gamma_{a} = \frac{2 - 12 \Gamma_{b}^{2} - 19 \Gamma_{b} \pm \sqrt{-116 \Gamma_{b} + 3 \Gamma_{b}^{2}(135 + 8 \Gamma_{b} (17 + 6 \Gamma_{b}))}}{2 + 24 \Gamma_{b}}
\end{gather}
This condition is only satisfied at $\Gamma_{a} \Gamma_{b} < 0$, $(\Gamma_{a}, \Gamma_{b}) = (1, 0)$, or $(\Gamma_{a}, \Gamma_{b}) = (0, 1)$, thus by the intermediate value theorem the determinant must be strictly non-negative.

To prove the roots are in the left half complex plane re-write the polynomial as
\begin{gather}
  z^{3} + a_{2} z^{2} + a_{1} z + a_{0} = 0\label{eq:hurwitz_poly}
\end{gather}
If $a_{0} = 0$, then since $a_{1}$ and $a_{2}$ are strictly positive we have a quadratic Hurwitz polynomial\cite{kuo1966}, and all the roots must have a strictly non-negative real component.
If $a_{0} \neq 0$, then if $a_{2} a_{1} > a_{0}$, \cref{eq:hurwitz_poly} is a Hurwitz polynomial.
This can be verified to be true by expanding $a_{2} a_{1}$ and considering each term separately.
Thus all roots $z$ must be strictly real non-positive.

As observed from \cref{fig:tf_eigs_dx_10} the stability restrictions imposed by the speed of light and $\mat{J}_{S}$ potentially can be orders of magnitude smaller than the stability requirements of the bulk ion motion.
This can be overly restrictive in situations where the fine scale dynamics associated with these faster speeds such as propagating Langmuir waves and plasma oscillations average out and have a negligible effect on the bulk plasma dynamics.
In addition to this, while computing an adaptive stable timestep in the advection dominated or source term dominated regimes is easy, it can be difficult to accurately choose a stable timestep in the transition regime.
Utilizing an A-stable or L-stable implicit time stepping scheme resolves both of these issues.

\section{Hybridizable discontinuous Galerkin methods}
\label{sec:hdg}


Despite being able to admit spatially discontinuous solutions, there is nothing at a theoretical mathematics level which prevents discontinuous Galerkin methods from being solved in an implicit fashion.
The primary practical challenge associated with solving classical DG discretizations with implicit time stepping is the large implicit global system which needs to be solved.
CG methods are able to mitigate some of the challenges associated with a large global implicit solve via static condensation\cite{guyan1965,glasser2004}.
HDG\cite{cockburn2009} generalizes this technique to DG methods by providing a secondary skeleton finite element space on element faces, to which element interior degrees of freedom are condensed to.
Since the degrees of freedom in each of these spaces scales with the surface area and volume respectively, for high enough order elements the condensed global skeleton space will end up with a smaller globally coupled system which needs to be solved.
Define the finite element spaces
\begin{gather}
  W_{h} := 
  \left\{
    w \in L^{2}(\Omega) : w|_{K} \in W(K) \forall K \in \Omega_{h}
  \right\}\\
  M_{h} := 
  \left\{
    m \in L^{2}(\partial \Omega_{h}) : m|_{e} \in M(e) \forall e \in \mathcal{E}_{h}, m|_{\partial \Omega_{D}} = 0
  \right\}
\end{gather}
where $W_{h}$ is the finite element space on elements and $M_{h}$ is the finite element space on faces.
Then consider a PDE system
\begin{gather}
  \vec{S}(t, \vec{q}, \vec{u}) + \diverge \mat{T}(t, \vec{q}, \vec{u}) = 0 \label{eq:hdg_prototype}\\
  \vec{V}(t, \vec{u}) = 0
\end{gather}
subject to the boundary conditions
\begin{gather}
  \vec{q}|_{\partial \Omega_{D}} = \vec{r}(t)\\
  \mat{T}(t, \vec{q}, \vec{u}) \cdot \hat{n}^{-}|_{\partial \Omega_{N}} = \vec{s}(t)
\end{gather}
where $\vec{q}$ may vary spatially while $\vec{u}$ are spatially invariant.
Treatment of any temporal derivatives for the purpose of the DG derivation is done separately via a method of lines\cite{schiesser1991} in the source terms operator $\vec{S}$.
For example, the two-fluid equations from \cref{sec:tf_model} would define
\begin{gather}
  \vec{S} =
  \begin{bmatrix}
    \partial_{t} \rho_{\alpha}\\
    \partial_{t} \vec{p}_{\alpha} - R_{c} \rho_{\alpha} (\vec{u}_{\beta} - \vec{u}_{\alpha}) - \frac{Z_{\alpha} L}{m_{\alpha} \delta_{p}} (\rho_{\alpha} \vec{E} + \vec{p}_{\alpha} \times \vec{B})\\
    \partial_{t} U_{\alpha} - Q_{c} (T_{\beta} - T_{\alpha})\\
    \partial_{t} \vec{E} + \frac{c_{0}^{2} L}{V_{A}^{2} \delta_{p}} \vec{j}\\
    \partial_{t} \vec{B}\\
    \frac{1}{c_{h}^{2}} \partial_{t} \psi + c_{p} \psi\\
    \frac{1}{c_{h}^{2}} \partial_{t} \theta + c_{p} \theta + \frac{c_{0}^{2} L}{V_{A}^{2} \delta_{p}}
  \end{bmatrix}
\end{gather}
The inclusion of $\vec{V}(t,\vec{u})$ is done to allow coupling to external multi-physics via boundary conditions. An example coupling could be an external circuit solver, where $\vec{u}$ would contain the circuit's degrees of freedom.

\begin{figure}[H]
  \centering{}
  \includegraphics[width=.33\textwidth]{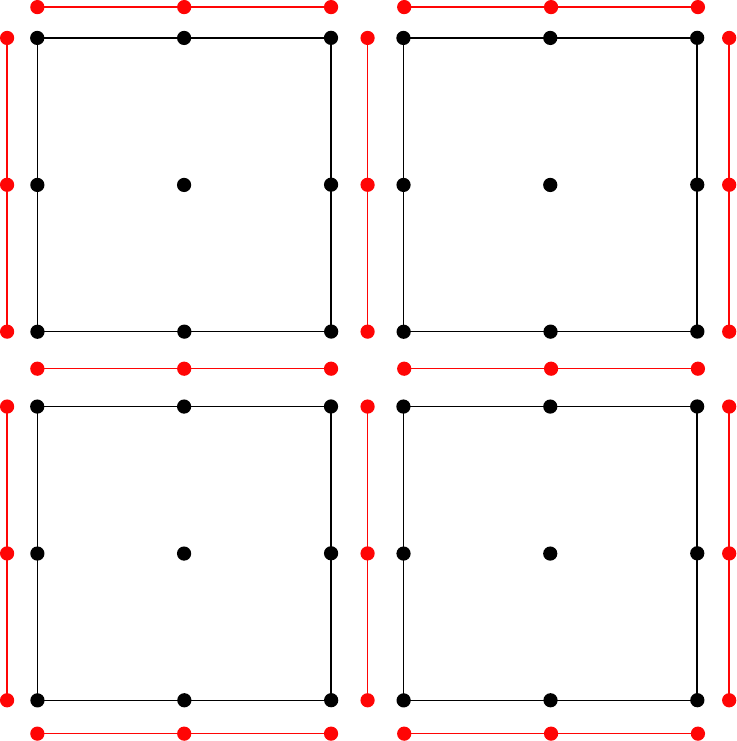}
  \caption{Example nodal degrees of freedom in HDG finite element spaces. Black nodes correspond to element interiors projected onto $W_{h}$ and red nodes correspond to the skeleton space projected onto $M_{h}$.}
  \label{fig:hdg_spaces}
\end{figure}

The HDG method then seeks to find an approximate solution $(\vec{q}_{h}, \vec{\lambda}_{h}, \vec{u})$ where each component of $\vec{q}_{h} \in W_{h}$ and $\vec{\lambda}_{h} \in M_{h}$ such that
\begin{gather}
  \vec{f} = \int_{\Omega_{h}} \vec{S}(\vec{q}_{h}, \vec{u}) v - \mat{T}(\vec{q}_{h}, \vec{u}) \cdot \grad v dV
  + \oint_{\partial \Omega_{h}} (\mat{T}(\vec{q}_{h}, \vec{\lambda}_{h}, \vec{u}) \cdot \hat{n})^{*} \mu dS = 0\\
  \vec{g} = \int_{\partial \Omega_{h}} \llbracket(\mat{T}(\vec{q}_{h}, \vec{\lambda}_{h}, \vec{u}) \cdot \hat{n})^{*} v\rrbracket dS - \int_{\partial \Omega_{n}} (\vec{s} - (\mat{T}(\vec{q}_{h}, \vec{\lambda}_{h}, \vec{u}) \cdot \hat{n})^{*}) v dS = 0\label{eq:conserve_flux}\\
  \vec{h} = \vec{V}(t, \vec{u}) = 0
\end{gather}
for all test functions $v \in W_{h}$ and $\mu \in M_{h}$, where $(\mat{T}(\vec{q}_{h}, \vec{\lambda}_{h}, \vec{u}) \cdot \hat{n})^{*}$ is the numerical flux.
That is to say, each scalar element of $\vec{q}_{h}$ is defined on the element interior finite element space $W_{h}$ while each scalar element of $\vec{\lambda}_{h}$ is defined on the skeleton facets finite element space $M_{h}$. An example nodal HDG discretization is shown in \cref{fig:hdg_spaces}.
No finite element space projection is required for the spatially invariant $\vec{u}$.
By choosing an appropriate numerical flux $(\mat{T} \cdot \hat{n})^{*}$ it becomes possible to decouple any direct dependency from one element to any other element; that is, the DOF of any element $\vec{q}_{i} \in \Omega_{i}$ are only directly coupled to the DOFs $\vec{q}_{i}$, $\vec{u}$, and $\vec{\lambda}_{j} \in e_{j}$ where $e_{j}$ are the faces of element $\Omega_{i}$.
Here we use the DG jump notation
\begin{gather}
  \llbracket
  q
  \rrbracket = q^{-} - q^{+}
\end{gather}
where $q^{-}$ is the interior value and $q^{+}$ is the exterior value, thus \cref{eq:conserve_flux} is a statement that the numerical flux across any interface must be conserved in a weak sense.

Since our goal is to solve coupled systems of non-linear PDEs, we will take the approach used by Nguyen and Peraire\cite{nguyen2012} for applying a Newton-Raphson iteration processes to HDG.
Take an initial guess for each unknown $\vec{q}_{h,0}$, $\vec{\lambda}_{h,0}$, and $\vec{u}_{0}$. Define a step size $\alpha \in (0,1]$ such that
\begin{gather}
  \Delta \vec{q}_{h} = \alpha (\vec{q}_{h,1} - \vec{q}_{h,0})\label{eq:ls_start}\\
  \Delta \vec{\lambda}_{h} = \alpha (\vec{\lambda}_{h,1} - \vec{\lambda}_{h,0})\\
  \Delta \vec{u} = \alpha (\vec{u}_{1} - \vec{u}_{0})\label{eq:ls_stop}
\end{gather}

The semi-discrete HDG problem can be written in matrix form as
\begin{gather}
  \begin{bmatrix}
    \mat{A} & \mat{B} & \mat{E}\\
    \mat{C} & \mat{D} & \mat{F}\\
    \mat{G} & \mat{H} & \mat{K}
  \end{bmatrix}
  \begin{bmatrix}
    \Delta \vec{q}_{h}\\
    \Delta \vec{\lambda}_{h}\\
    \Delta \vec{u}
  \end{bmatrix} =
  \begin{bmatrix}
    -\vec{f}_{0}\\
    -\vec{g}_{0}\\
    -\vec{h}_{0}
  \end{bmatrix}\\
  \begin{aligned}
    \mat{A} &= \partial_{\vec{q}} \vec{f}_{0} & \mat{B} &= \partial_{\vec{\lambda}} \vec{f}_{0} & \mat{E} &= \partial_{\vec{u}} \vec{f}_{0}\\
    \mat{C} &= \partial_{\vec{q}} \vec{g}_{0} & \mat{D} &= \partial_{\vec{\lambda}} \vec{g}_{0} & \mat{F} &= \partial_{\vec{u}} \vec{g}_{0}\\
    \mat{G} &= \partial_{\vec{q}} \vec{h}_{0} & \mat{H} &= \partial_{\vec{\lambda}} \vec{h}_{0} & \mat{K} &= \partial_{\vec{u}} \vec{h}_{0}
  \end{aligned}
\end{gather}
where the 0 subscript is used to denote evaluations at the initial guess.
Since $\mat{A}$ is guaranteed to be block-diagonal, it is relatively cheap to invert. Thus it becomes feasible to construct the Schur complement system
\begin{gather}
  \begin{bmatrix}
    \mat{D} - \mat{C} \mat{A}^{-1} \mat{B} & \mat{F} - \mat{C} \mat{A}^{-1} \mat{E}\\
    \mat{H} - \mat{G} \mat{A}^{-1} \mat{B} & \mat{K} - \mat{G} \mat{A}^{-1} \mat{E}
  \end{bmatrix}
  \begin{bmatrix}
    \Delta \vec{\lambda}_{h}\\
    \Delta \vec{u}
  \end{bmatrix}=
  \begin{bmatrix}
    -\vec{g} + \mat{C} \mat{A}^{-1} \vec{f}\\
    -\vec{h} + \mat{G} \mat{A}^{-1} \vec{f}
  \end{bmatrix}
\end{gather}
The blocks of the Schur complement have the following properties:
\begin{enumerate}
\item $\mat{D} - \mat{C} \mat{A}^{-1} \mat{B}$ has a well-defined sparsity pattern proportional to the number of faces shared between two connected elements\cite{peraire2008,kirby2012}.
\item the remaining 3 blocks are in general dense, however have at least 1 dimension which has the same length as $\vec{u}$, which is assumed to be very small.
\end{enumerate}

After solving the Schur system, the element interior step direction can be computed using
\begin{gather}
  \Delta \vec{q}_{h} = -\mat{A}^{-1} \vec{f} - \mat{A}^{-1} \mat{B} \Delta \vec{\lambda}_{h} - \mat{A}^{-1} \mat{E} \Delta \vec{u}
\end{gather}
Given step directions $\Delta \vec{q}_{h}$, $\Delta \lambda_{h}$, and $\Delta u$ we can then apply a linesearch\cite{dennis1996} to compute $\alpha$, improving global convergence properties of the base Newton-Raphson method. Note that fixing $\alpha = 1$ gives the traditional Newton Raphson method. Finally, we can update our guess by solving \cref{eq:ls_start}-\cref{eq:ls_stop} for $\vec{q}_{h,1}$, $\vec{\lambda}_{h,1}$, and $\vec{u}_{1}$.
This process is iterated until a sufficiently accurate solution is found, judged using a combination of $\|f_{0}\| + \|g_{0}\| + \|h_{0}\| \approx 0$ and $\|\Delta \vec{q}_{h} \| + \|\Delta \vec{\lambda}_{h}\| + \|\Delta \vec{u}\| \approx 0$.

\subsection{Problem decomposition and implementation details}


To simplify and efficiently recovery the element interior solutions, define the intermediate quantities
\begin{gather}
  \mat{\Gamma} = \mat{A}^{-1} \mat{B}\\
  \mat{\Omega} = \mat{A}^{-1} \mat{E}\\
  \vec{\xi} = \mat{A}^{-1} \vec{f}
\end{gather}
$\mat{\Gamma}$ is block-sparse, with only the blocks $\mat{\Gamma}_{ef}$ being non-zero, where $e$ is a given element's row and $f$ are columns associated with faces neighboring face $e$.
In general $\mat{\Omega}$ is dense, however the width is of length $\vec{u}$, which is assumed to be small.
This allows the element interior solution to be recovered as
\begin{gather}
  \Delta \vec{q}_{h} = - \vec{\xi} - \mat{\Gamma} \Delta \vec{\lambda}_{h} - \mat{\Omega} \Delta u
\end{gather}
and avoids recomputing $\mat{A}^{-1}$. For high order elements and large coupled systems this calculation, while still block diagonal, is non-trivial and involves inverting block dense matrix with hundreds of rows.
Each of these matrices can be computed in an element-by-element fashion, meaning that $\mat{A}$, $\mat{C}$, and $\mat{G}$ never need to be formed in their entirety, but rather only the portions associated with a given element and its faces are required at any one time.
The full non-linear HDG solve can then be implemented along with a line search\cite{dennis1996} as:
\begin{breakablealgorithm}
  \begin{algorithmic}
    \Function{Nonlinear HDG}{$\vec{q}$, $\vec{\lambda}$, $\vec{u}$}
      \While{solution not converged}
        \State{$\mat{S} \gets 0$, $\vec{b} \gets 0$, where $\mat{S}$ and $\vec{b}$ constitute the storage space for the Schur complement system $\mat{S} \vec{\Delta} = \vec{b}$.}
        \State{Add $\mat{K}(\vec{f}, \vec{\lambda}, \vec{u})$ to the appropriate rows/cols of $\mat{S}$}
        \State{Subtract $\vec{h}(\vec{q}, \vec{\lambda}, \vec{u})$ from the appropriate rows of $\vec{b}$}
        \For{each element $i$}
          \State{$\mat{A}^{-1} \gets \mat{A}^{-1}(\vec{q}_{i}, \vec{u})$}
          \For{each face $j$ connected to element $i$}
            \State{$\mat{\Gamma}_{ij} \gets \mat{A}^{-1} \mat{B}(\vec{q}_{i}, \vec{\lambda}_{j}, \vec{u})$}
            \State{$\mat{C}_{j} \gets \mat{C}(\vec{q}_{i}, \vec{\lambda}_{j}, \vec{u})$}
            \State{Add $\mat{D}(\vec{q}_{i}, \vec{\lambda}_{j}, \vec{u})$, $\mat{F}(\vec{q}_{i}, \vec{u}_{j}, \vec{u})$, and $\mat{H}(\vec{q}_{i}, \vec{\lambda}_{j}, \vec{u})$, to the appropriate rows/cols of $\mat{S}$}
            \State{Subtract $\vec{g}(\vec{q}_{i}, \vec{\lambda}_{j}, \vec{u})$ from the appropriate rows of $\vec{b}$}
          \EndFor{}
          \State{$\vec{f} \gets \vec{f}(\vec{q}_{i}, \vec{\lambda}, \vec{u})$, $\mat{G} \gets \vec{G}(\vec{q}_{i}, \vec{\lambda}, \vec{u})$, $\mat{E} \gets \mat{E}(\vec{q}_{i}, \vec{\lambda}, \vec{u})$}
          \State{$\vec{\xi}_{i} \gets -\mat{A}^{-1} \vec{f}$}
          \State{$\mat{\Omega}_{i} \gets \mat{A}^{-1} \mat{E}$}
          \State{Subtract $\mat{G} \vec{\xi}_{i}$ from the appropriate rows of $\vec{b}$}
          \State{Subtract $\mat{G} \mat{\Omega}_{i}$ from the appropriate rows/cols of $\mat{S}$}
          \For{each face $j$ connected to element $i$}
            \State{Subtract $\mat{G} \mat{\Gamma}_{ij}$ from the appropriate rows/cols of $\mat{S}$}
            \For{each face $k$ connected to element $i$}
              \State{Subtract $\mat{C}_{j} \mat{\Gamma}_{ik}$ from the appropriate rows/cols of $\mat{S}$}
            \EndFor{}
          \EndFor{}
        \EndFor{}
        \State{Solve the linearized system
$
            \mat{S}
            \begin{bmatrix}
              \Delta \vec{\lambda}\\
              \Delta \vec{u}
            \end{bmatrix} = \vec{b}
          $
        }
        \For{each element $i$}
          \State{$\Delta q_{i} \gets \vec{\xi}_{i} - \mat{\Omega}_{i} \Delta \vec{u}$}
          \For{each face $j$ connected to element $i$}
            \State{$\Delta q_{i} \gets \Delta q_{i} - \mat{\Gamma}_{ij} \Delta \vec{\lambda}_{j}$}
          \EndFor{}
        \EndFor{}
        \State{$\alpha \gets \text{linesearch}(\vec{q}, \vec{\lambda}, \vec{u}, \Delta \vec{q}, \Delta \vec{\lambda}, \Delta \vec{u})$}
        \State{$\vec{q} \gets \vec{q} + \alpha \Delta \vec{q}$, $\vec{\lambda} \gets \vec{\lambda} + \alpha \Delta \vec{\lambda}$, $\vec{u} \gets \vec{u} + \alpha \Delta \vec{u}$}
      \EndWhile{}
      \State{\Return{$\vec{q}$, $\vec{\lambda}$, $\vec{u}$}}
    \EndFunction{}
  \end{algorithmic}
\end{breakablealgorithm}

~\\
For a large PDE system such as the two-fluid plasma model, deriving correct analytical Jacobians by hand is practically impossible due to the large number of non-zero entries.
For the two-fluid plasma system this requires at least 1047 unique non-zero entries, before taking into account boundary conditions which may double or triple the number of analytical Jacobian non-zeroes.
To solve this problem the Jacobian along with the rest of the HDG code is generated using SymPy\cite{sympy}, a symbolic computer algebra system.
The tool is capable of automatically detecting and eliminating unused degrees of freedom and handles multiple subdomains/boundary conditions.
The generated HDG code is coupled with a diagonally implicit Runge-Kutta\cite{hairer1980} temporal solver, and parallelized on CPUs using MPI.

\subsection{Face-oriented numerical flux}


A key requirement for the HDG numerical flux is that from the perspective of any element the numerical flux must be a function of only the DOFs inside the element and the DOF on the neighboring face.
Each face is assigned a local normal, tangent, and binormal. As a result the face normal will point outwards for one element, while for the adjacent neighboring element it points inwards.
The non-uniqueness of the tangent and binormal directions is not an issue so long as they are chosen such that $\hat{t} = \hat{n} \times \hat{b}$.
Vector and tensor quantities are appropriately rotated to the face coordinate system for computing numerical fluxes.
\begin{figure}[H]
  \centering{}
  \includegraphics{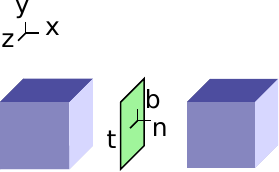}
  \caption{Global and face-oriented coordinate systems}
\end{figure}
Computing the numerical flux becomes
\begin{gather}
  (\mat{T} \hat{n})^{*} = \pm \vec{T}_{n}^{*}
\end{gather}
Depending on the given system this could potentially lead to a reduction of face DOF\cite{bui-thanh2016}.
For example, take the linear diffusion equation
\begin{gather}
  \partial_{t} \phi - \diverge \vec{\gamma} = 0\\
  \vec{\gamma} = \grad \phi
\end{gather}
The numerical fluxes chosen are the class of central fluxes of the form
\begin{gather}
  (\mat{T} \hat{n})^{*} =
  \mat{T}_{\lambda} \hat{n} + \tau (\vec{q}^{-} - \vec{\lambda})
\end{gather}
For the linear diffusion equation this expands to
\begin{gather}
  (\mat{T} \hat{n})^{*} =
  \begin{bmatrix}
    \lambda_{\vec{\gamma}} \cdot \hat{n}^{-} + \tau (q^{-}_{\phi} - \lambda_{\phi})\\
    -\lambda_{\phi} \mat{I} \hat{n}^{-} + \tau (q^{-}_{\vec{\gamma}} - \lambda_{\vec{\gamma}})
  \end{bmatrix}
\end{gather}
where $\tau > 0$ is a user chosen constant. Note that in practice so long as $\tau \approx 1$ the method will converge at the same rate.
In a global coordinate system this expands to
\begin{gather}
  (\mat{T} \hat{n})^{*} =
  \begin{bmatrix}
    \lambda_{\gamma_{x}} \hat{n}^{-}_{x} + \lambda_{\gamma_{y}} \hat{n}^{-}_{y} + \lambda_{\gamma_{z}} \hat{n}^{-}_{z} + \tau (q^{-}_{\phi} - \lambda_{\phi})\\
    \begin{bmatrix}
      -\lambda_{\phi} \hat{n}^{-}_{x} + \tau (q^{-}_{\gamma_{x}} - \lambda_{\gamma_{x}})\\
      -\lambda_{\phi} \hat{n}^{-}_{y} + \tau (q^{-}_{\gamma_{y}} - \lambda_{\gamma_{y}})\\
      -\lambda_{\phi} \hat{n}^{-}_{z} + \tau (q^{-}_{\gamma_{z}} - \lambda_{\gamma_{z}})
    \end{bmatrix}
  \end{bmatrix}
\end{gather}
which requires four element scalar fields and four face scalar fields.
In a rotated frame this can be simplified to
\begin{gather}
  (\mat{T} \hat{n})^{*} =
  \begin{bmatrix}
    \pm \lambda_{\gamma_{n}} + \tau (q^{-}_{\phi} - \lambda_{\phi})\\
    \begin{bmatrix}
      \mp \lambda_{\phi} + \tau (q^{-}_{\gamma_{n}} - \lambda_{\gamma_{n}})\\
      0\\
      0
    \end{bmatrix}
  \end{bmatrix}
\end{gather}
eliminating two face scalar fields.
For the two-fluid plasma model with either mixed parabolic-hyperbolic cleaning or purely parabolic cleaning this rotation scheme requires a total of 42 element and 34 face scalar fields, eliminating eight face scalar fields from the non-rotated scheme.

\section{Validation of convergence rates for various linear systems}
\label{sec:convergence}

\begin{figure}[H]
  \centering{}
  \includegraphics[width=.5\textwidth]{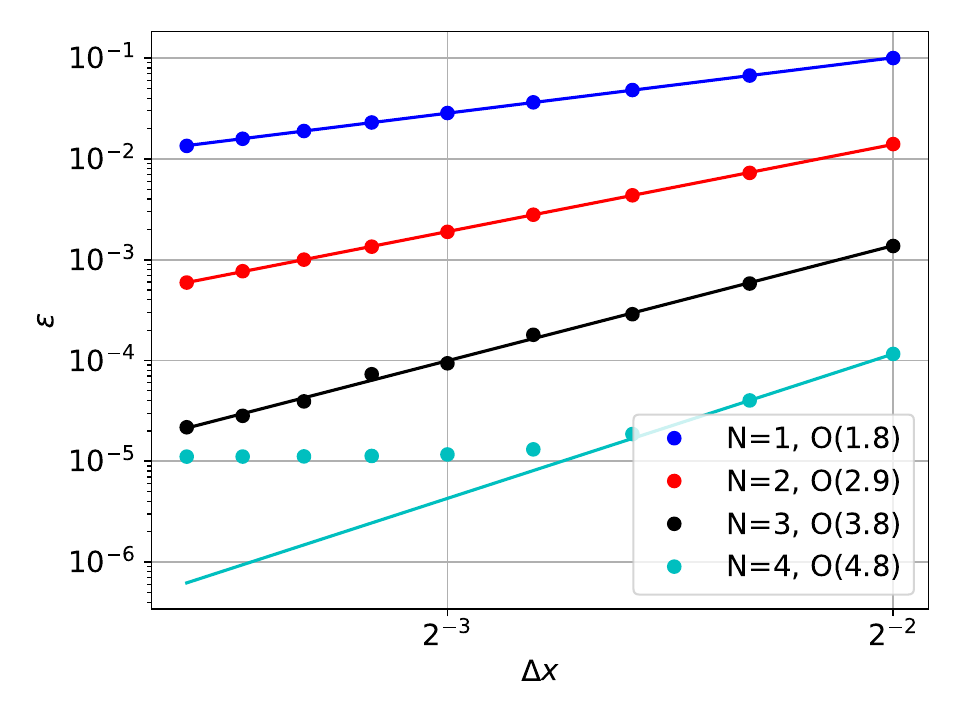}
  \caption{L-2 error norms for the linear advection equation converge at the expected optimal $O(N+1)$ rates for degree $N$ basis. The observed $10^{-5}$ minimum error is from the temporal discretization error.}
  \label{fig:advection_convergence}
\end{figure}

Convergence is evaluated for several model PDE systems including the linear advection, diffusion equation, and wave equation. In all cases periodic boundary conditions are applied on a 3D cube domain $\vec{r} \in [0,1]^{3}$ and $t \in [0, 0.1]$ with a fixed $\Delta t = 5 \times 10^{-3}$.
The error is compared against the analytical solution $a_{a}$ using an $8^{3}$ tensor product Gauss-Lobatto quadrature to compute the L-2 error integral
\begin{gather}
  \epsilon = \sqrt{\int (q_{h} - q_{a})^{2} dV}
\end{gather}
The linear advection system is defined by
\begin{gather}
  \partial_{t} q + \diverge (\vec{a} q) = 0\\
  \vec{a}^{T} =
  \begin{bmatrix}
    1 & \frac{1}{2} & 2
  \end{bmatrix}\\
  q_{a}(x,y,z,t) = \cos(2\pi (x - a_{x} t)) \cos(2\pi (y - a_{y} t)) \cos(2\pi (z - a_{z} t))
\end{gather}
\Cref{fig:advection_convergence} demonstrates that the linear advection system converges at the expected optimal $O(N+1)$ rates.
Note that at small $\Delta x$ and high basis order the simulation reaches a saturation floor of $10^{-5}$ from the temporal discretization error. This can be reduced by choosing a smaller $\Delta t$ or more accurate temporal discretization method.

\begin{figure}[H]
  \centering{}
  \begin{subfigure}{.4\textwidth}
    \includegraphics[width=\textwidth]{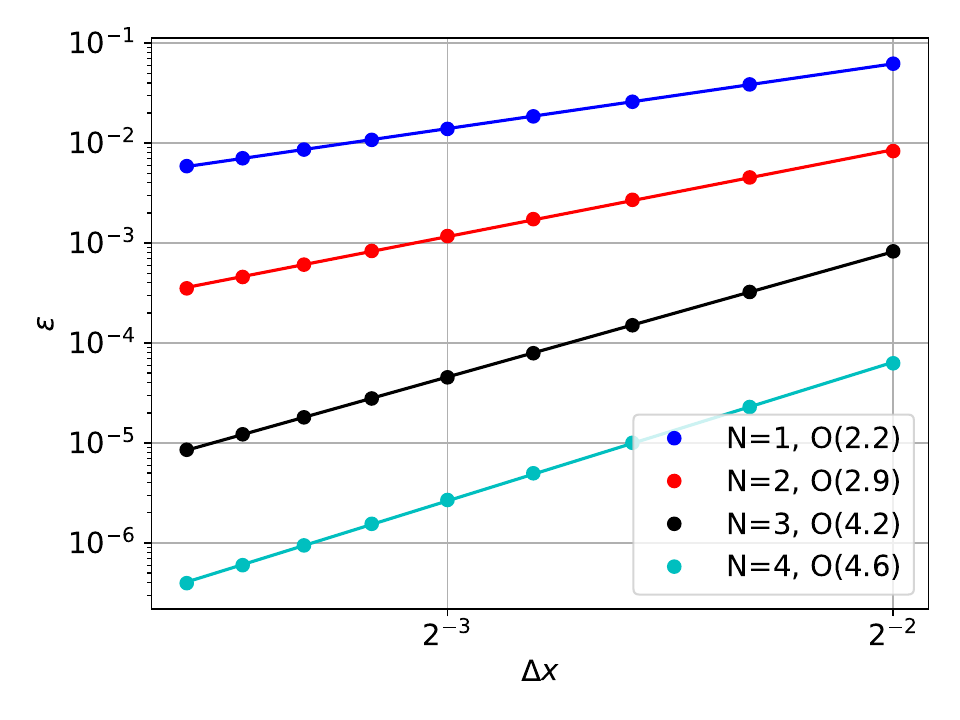}
    \caption{$q$}
  \end{subfigure}
  \begin{subfigure}{.4\textwidth}
    \includegraphics[width=\textwidth]{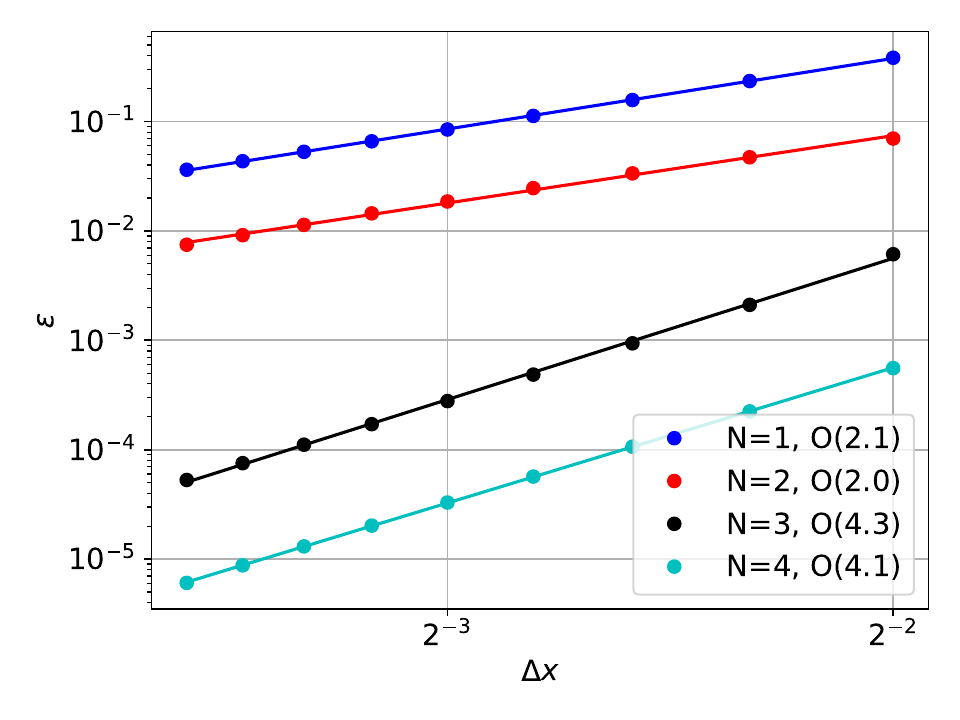}
    \caption{$u_{x}$}
  \end{subfigure}
  \caption{L-2 error norms for the linear diffusion equation. Note that $q$ converges at $O(N+1)$, while $\vec{u}$ has an odd-even oscillation between $O(N+1)$ and $O(N)$ convergence rates. $u_{y}$ and $u_{z}$ have near identical plots to $u_{x}$ and are omitted for brevity.}
  \label{fig:diffusion_convergence}
\end{figure}

The linear diffusion system is defined by
\begin{gather}
  \partial_{t} q - \diverge (k \grad q) = 0
\end{gather}
This can be decomposed into a coupled system of first order equations to fit the form of \cref{eq:hdg_prototype} as
\begin{gather}
  \partial_{t} q - \diverge (k \vec{u}) = 0\\
  \vec{u} - \grad q = 0\\
  k = 10^{-2}\\
  q_{a}(x,y,z,t) = \exp
  \left(
    -3 (2\pi)^{2} k t
  \right) \sin(2\pi x) \sin(2\pi x) \sin(2\pi z)
\end{gather}
where initial conditions are given by $q(x,y,z) = q_{a}(x,y,z,t=0)$.
\Cref{fig:diffusion_convergence} demonstrates optimal $O(N+1)$ convergence rates for $q$, however $\vec{u}$ converges with a mix of $O(N+1)$ and $O(N)$ rates depending if $N$ is odd or even.
This oscillating convergence rate is a known behavior of central fluxes\cite{hesthaven2008,cockburn1997}. Implementing a bidirectional upwinding flux for handling second order spatial derivatives has been observed to resolve the sub-optimal convergence rates for even $N$\cite{cockburn2009}.

\begin{figure}[H]
  \centering{}
    \begin{subfigure}{.32\textwidth}
    \includegraphics[width=\textwidth]{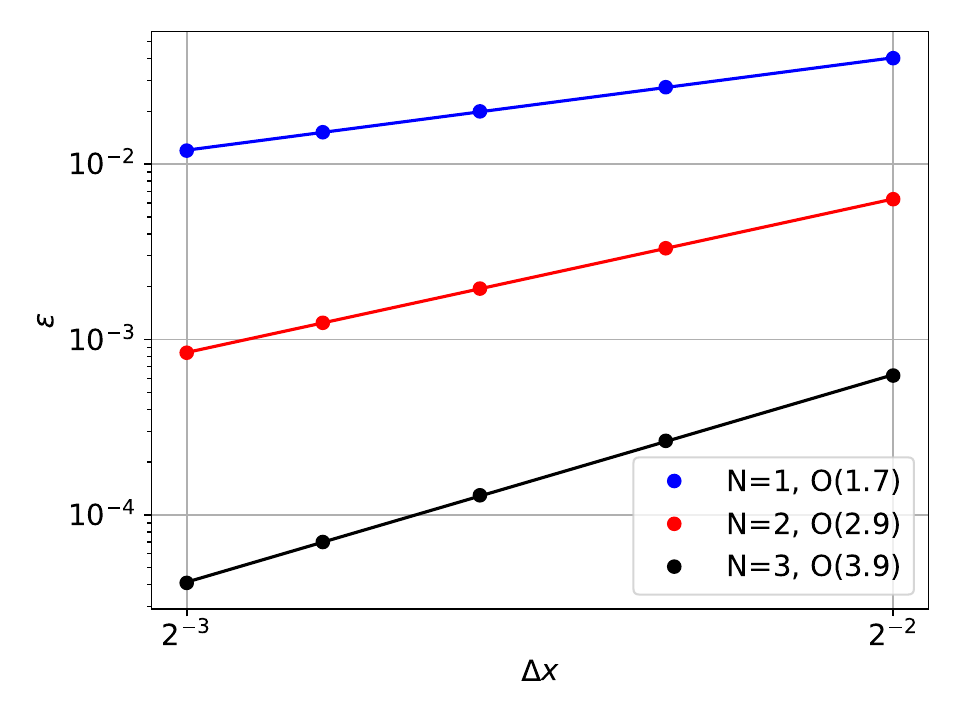}
    \caption{$q$}
  \end{subfigure}
  \begin{subfigure}{.32\textwidth}
    \includegraphics[width=\textwidth]{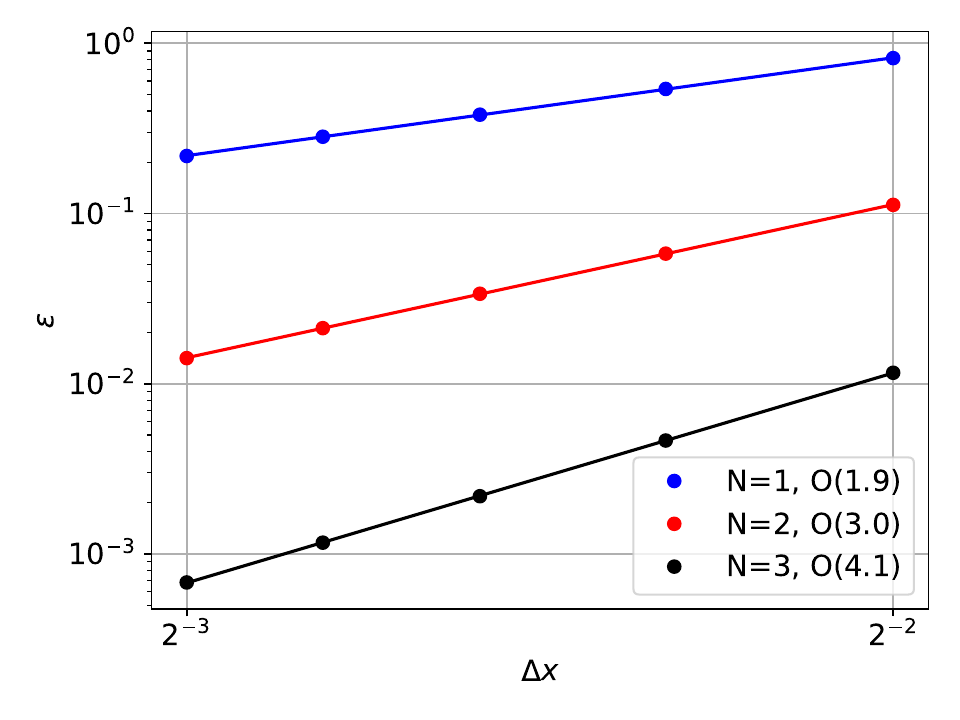}
    \caption{$v$}
  \end{subfigure}
  \begin{subfigure}{.32\textwidth}
    \includegraphics[width=\textwidth]{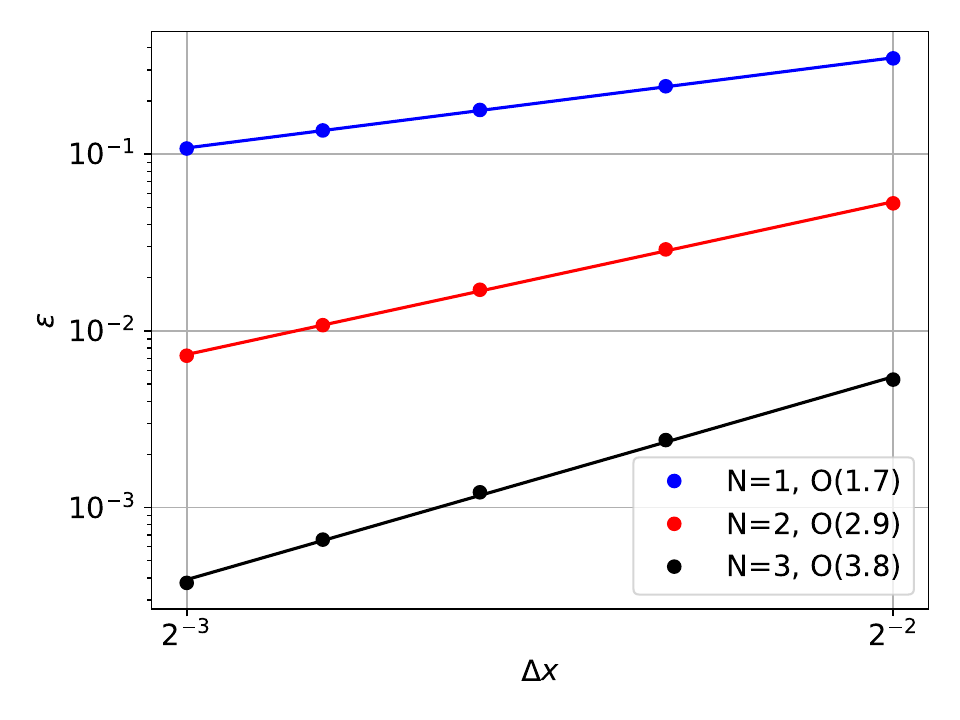}
    \caption{$u_{x}$}
  \end{subfigure}
  \caption{L-2 error norms for the linear wave equation. All variables converge at the optimal $O(N+1)$ rate. Plots for $u_{y}$ and $u_{z}$ are nearly identical to $u_{x}$ and omitted.}
  \label{fig:wave_convergence}
\end{figure}

The linear wave equation is defined by
\begin{gather}
  \partial_{t} q - v = 0\label{eq:wave1}\\
  \partial_{t} \vec{u} - \grad v = 0\label{eq:wave2}\\
  \partial_{t} v - \diverge(c^{2} \vec{u}) = 0\label{eq:wave3}\\
  c = 1\\
  q_{a}(x,y,z,t) = \cos
  \left(
    2\pi t \sqrt{3 c^{2}}
  \right)
  \sin(2\pi x) \sin(2\pi y) \sin(2\pi z)
\end{gather}
We can recover the traditional wave equation by solving \cref{eq:wave1} for $v$ and substituting into \cref{eq:wave2} and \cref{eq:wave3}:
\begin{gather}
  \partial_{t} \vec{u} - \grad \partial_{t} q = 0\nonumber\\
  \begin{aligned}
    \partial_{t} (\vec{u} - \grad q) &= 0 & \rightarrow && \vec{u} &= \grad q
  \end{aligned}
  \nonumber\\
  \partial_{tt} q - \diverge
  \left(
    c^{2} \grad q
  \right) = 0
\end{gather}
We choose to not solve the wave equation in this form because of the presence of second order time and spatial derivatives, which our implementation cannot directly handle.
\Cref{fig:wave_convergence} demonstrates that all variables present in the simulation converge at the optimal $O(N+1)$ rate. Unlike the diffusion equation the use of central fluxes does not negatively impact convergence rates despite using the same splitting method for the Laplacian operator\cite{cockburn2016}.

\section{Application of HDG to the two-fluid plasma system}
\label{sec:plasma_tests}

\begin{figure}[H]
  \centering{}
  \begin{subfigure}{.49\textwidth}
    \includegraphics[width=\textwidth]{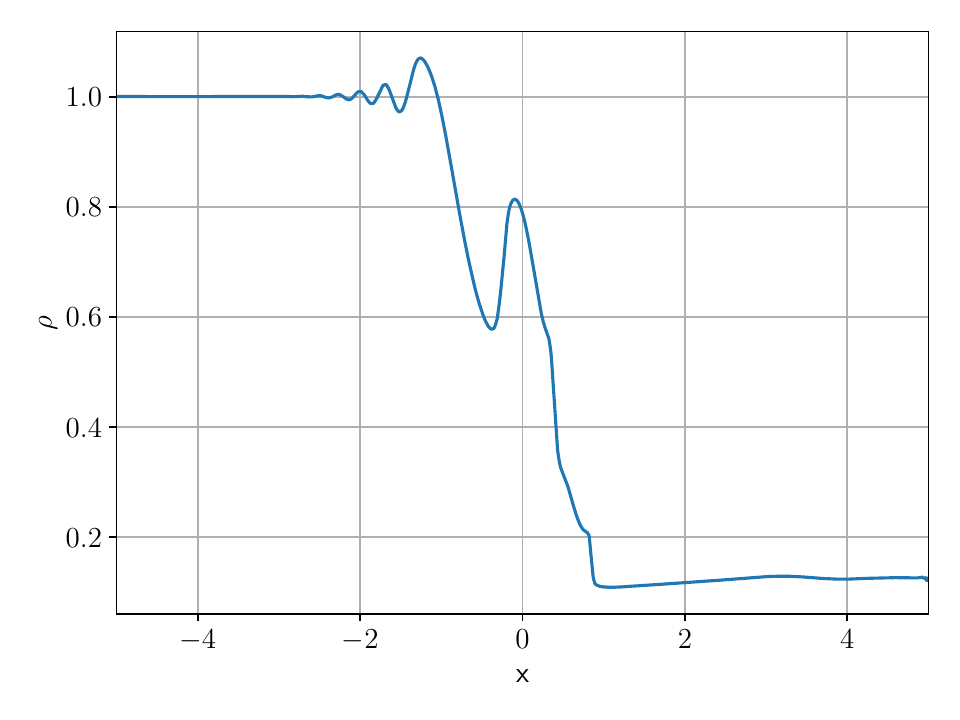}
    \caption{$\rho$}
  \end{subfigure}
  \begin{subfigure}{.49\textwidth}
    \includegraphics[width=\textwidth]{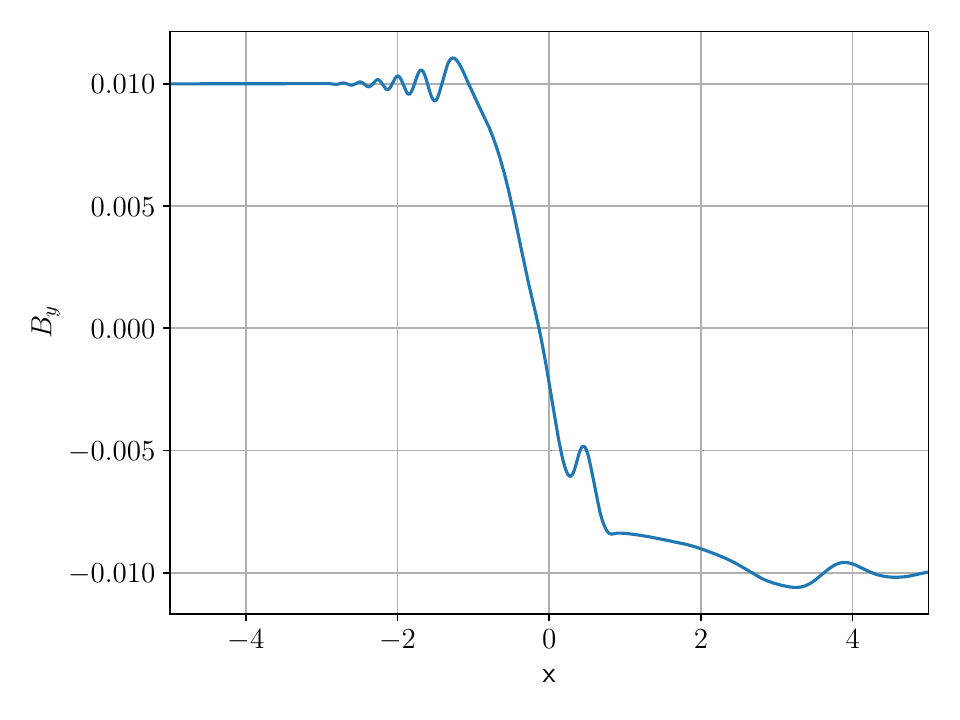}
    \caption{$B_{y}$}
  \end{subfigure}
  \begin{subfigure}{.49\textwidth}
    \includegraphics[width=\textwidth]{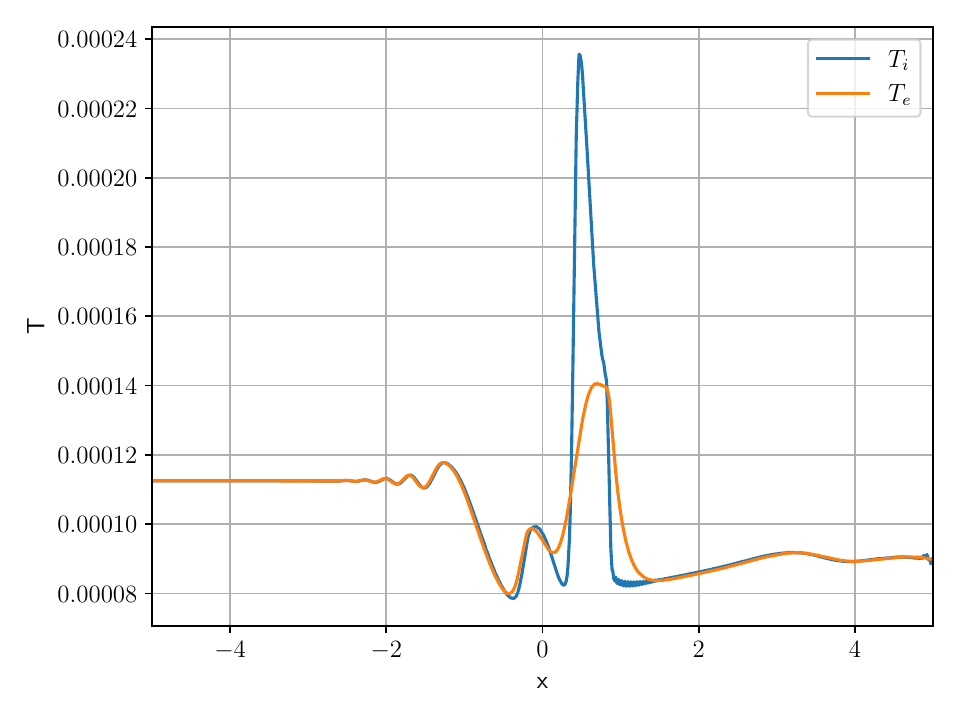}
    \caption{$T$}
  \end{subfigure}
  \begin{subfigure}{.49\textwidth}
    \includegraphics[width=\textwidth]{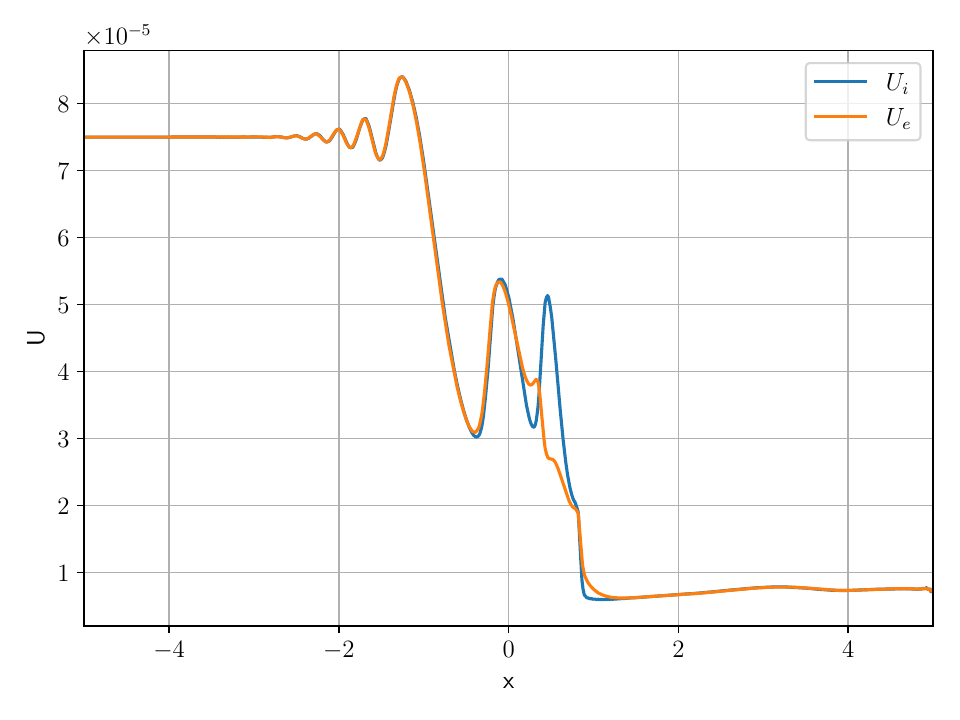}
    \caption{$U$}
  \end{subfigure}
  \caption{Two-fluid Brio-Wu shocktube properties at $t=100$. Similar features to those found in Shumlak and Loverich\cite{shumlak2003} are found, which are notably distinct from the classical MHD results.}
  \label{fig:brio_wu}
\end{figure}

The two-fluid plasma system was used to test the behavior of the HDG code on large coupled systems.
A quasi-1D domain of $256\times 1\times 1$ parabolic elements with $x \in [-5, 5]$ is initialized with a two-fluid extension of the magnetized Brio-Wu shock tube\cite{shumlak2003, loverich2011} such that
\begin{gather}
  \begin{aligned}
  \begin{bmatrix}
    n_{i}\\
    U_{i}\\
    n_{e}\\
    U_{e}\\
    B_{x}\\
    B_{y}
  \end{bmatrix}_{\text{left}} &=
                  \begin{bmatrix}
                    1\\
                    7.5 \times 10^{-5}\\
                    1\\
                    7.5 \times 10^{-5}\\
                    7.5 \times 10^{-3}\\
                    10^{-2}
                  \end{bmatrix} &
  \begin{bmatrix}
    n_{i}\\
    U_{i}\\
    n_{e}\\
    U_{e}\\
    B_{x}\\
    B_{y}
  \end{bmatrix}_{\text{right}} &=
                  \begin{bmatrix}
                    0.125\\
                    7.5 \times 10^{-5}\\
                    0.125\\
                    7.5 \times 10^{-6}\\
                    7.5 \times 10^{-3}\\
                    -10^{-2}
                  \end{bmatrix}
  \end{aligned}\\
  \begin{aligned}
    \frac{\delta_{p}}{L} &= 1 & \gamma &= \frac{5}{3} & m_{i} &= 1 & m_{e} &= \frac{1}{1836} & \frac{c_{0}}{V_{A}} &= c_{h} = c_{p} = 1 & Z_{i} &= -Z_{e} = 10
  \end{aligned}\\
  \mu_{c} = \kappa_{c} = Q_{c} = R_{c} = 10^{-5}
\end{gather}
The is simulated from $t \in [0,100]$.
\Cref{fig:brio_wu} shows various plasma properties at the end of the simulation.
The plotted fields show similar features to the work of Shumlak and Loverich\cite{shumlak2003}, while being distinct from results produced using standard MHD.
Some examples of these differences is the presence of a Whistler wave propagating to the left around $x=1.75$, and the solution producing $T_{i} \neq T_{e}$.
To validate the correct propagation speed of various important two-fluid waves, a second simulation with small jumps is initialized to ensure linear behavior.
The domain is changed to $1024\times 1 \times 1$ linear elements with $x \in [-\pi, \pi]$, and
\begin{gather}
  \begin{aligned}
  \begin{bmatrix}
    n_{i}\\
    U_{i}\\
    n_{e}\\
    U_{e}\\
    B_{x}\\
    B_{y}
  \end{bmatrix}_{\text{left}} &=
                  \begin{bmatrix}
                    10.0005\\
                    7.500375 \times 10^{-5}\\
                    10.0005\\
                    7.500375 \times 10^{-5}\\
                    7.5 \times 10^{-2}\\
                    5 \times 10^{-5}
                  \end{bmatrix} &
  \begin{bmatrix}
    n_{i}\\
    U_{i}\\
    n_{e}\\
    U_{e}\\
    B_{x}\\
    B_{y}
  \end{bmatrix}_{\text{right}} &=
                  \begin{bmatrix}
                    9.9995\\
                    7.499625 \times 10^{-5}\\
                    9.9995\\
                    7.499625 \times 10^{-5}\\
                    7.5 \times 10^{-2}\\
                    -5 \times 10^{-5}
                  \end{bmatrix}
  \end{aligned}\\
  \begin{aligned}
    \frac{\delta_{p}}{L} &= 1 & \gamma &= \frac{5}{3} & m_{i} &= 1 & 1836 m_{e} &= 1 & \frac{c_{0}}{V_{A}} &= c_{h} = c_{p} = 1 & Z_{i} &= -Z_{e} = 1
  \end{aligned}\\
  \mu_{c} = \kappa_{c} = Q_{c} = R_{c} = 10^{-5}
\end{gather}
A mixed 2D spatial/temporal Fourier analysis is then performed on $E_{y}$\cite{shumlak2003}.
This is plotted against the theoretical L and R mode wave propagation speeds, described by the dispersion relation
\begin{gather}
  \frac{c_{0}^{2} k^{2}}{V_{A}^{2} \omega^{2}} = 1 - \frac{\omega_{pe}^{2}}{\omega (\omega \pm \omega_{ce})}\\
  \begin{aligned}
    \omega_{pe} &= (\omega_{p} \tau) \sqrt{\frac{n_{e} Z_{e}^{2}}{m_{e}}}
    &
      \omega_{ce} &= \frac{L}{\delta_{p}} \frac{|Z_{e}| B}{m_{e}}
  \end{aligned}
\end{gather}
\begin{figure}[H]
  \centering{}
  \includegraphics[width=.6\textwidth]{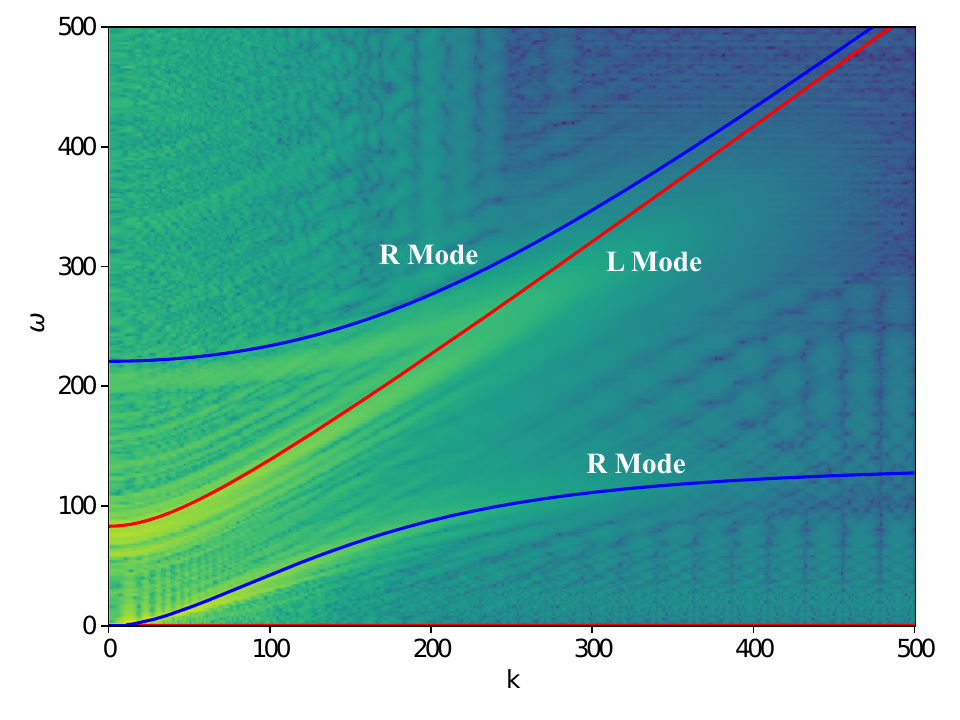}
  \caption{Fourier analysis of $E_{y}$ matches the circularly polarized two-fluid plasma dispersion relation.}
  \label{fig:dispersion}
\end{figure}

\Cref{fig:dispersion} demonstrates that the Fourier analysis has three dominant modes corresponding to the three branches of the L and R modes.
This test demonstrates the viability of using HDG to solve large-scale coupled plasma differential equations and the validity of the implemented tools to ensure correctness and simplify implementing such solvers.

\begin{figure}[H]
  \centering{}
  \includegraphics[width=.6\textwidth]{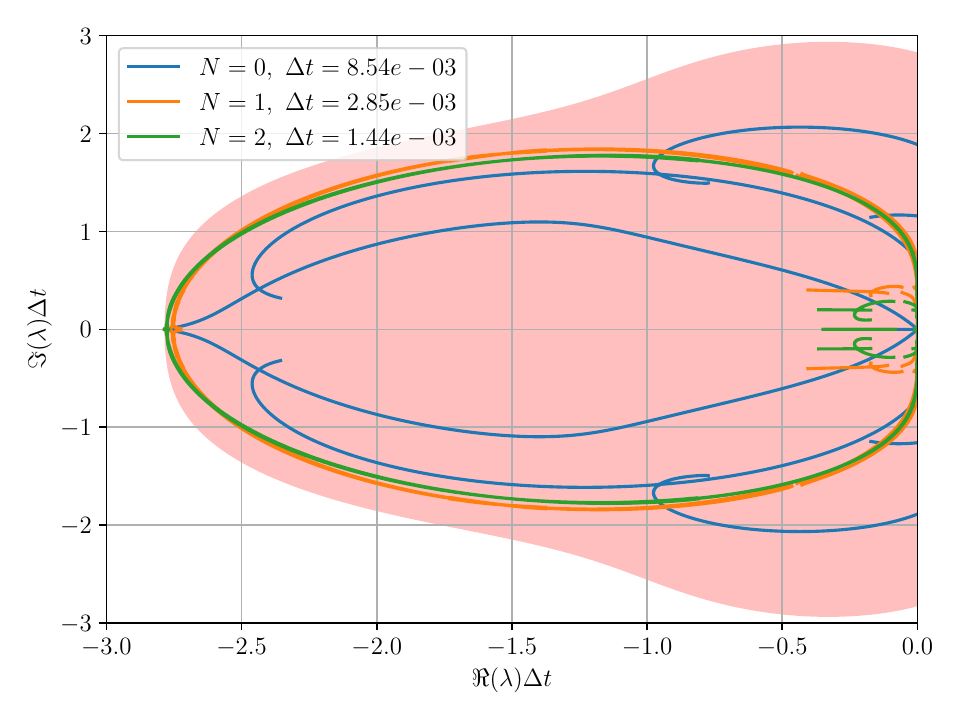}
  \caption{Stability conditions for a two-fluid explicit RK4 DG solver for the above specified conditions. The red highlighted region is the region of absolute stability for RK4 and the curves define the stability eigencurves of the ideal two-fluid system generated using the method from \cref{sec:tf_model}. Implicit simulations were run with $\Delta t = 6.28 \times 10^{-3}$ to achieve sufficient temporal resolution for the Fourier dispersion analysis in \cref{fig:dispersion}. A degree 0 basis could be run explicitly with this timestep; however, an explicit solution with a high degree basis is unstable using this timestep.}
  \label{fig:lr_num_stability}
\end{figure}

\Cref{fig:lr_num_stability} shows the numerical stability conditions for an equivalent explicit explicit RK4. For analysis purposes the implicit scheme used $\Delta t = 6.28 \times 10^{-3}$. While an explicit RK4 DG solver can run a $N=0$ basis with this timestep, higher degree basis require a smaller timestep.
Define the timestep ratio gain as the implicit timestep $\Delta t$ used divided by the maximum stable explicit timestep.
In this case there is a timestep ratio gain of $2.2$ for $N=1$ and $4.4$ for $N=2$ over an explicit RKDG method.
The primary reason these ratios are low is that the original test conditions were setup such that it is feasible to run and benchmark an explicit RKDG code.
Nevertheless, it is sufficient for demonstrating that the HDG method can continue to operate accurately in a regime where explicit time stepping methods become unstable.
In other test problems we have observed stable and accurate runs using timesteps greater than $60$ times larger than allowable for the explicit RK4 DG method.

\section{Conclusion}


This work increases the viability of using the two-fluid plasma model, which while more physically accurate than single fluid models is numerically stiff\cite{srinivasan2011}.
We perform a detailed numerical stability analysis of solving the ideal two-fluid system using classical RKDG. The numerically stiffness of the two-fluid plasma model results from needing to temporally resolve the speed of light, $c_{0}/V_{A}$, and the eigenvalues of the source terms, $\mat{J}_{S}$, when bulk plasma motion is dictated by timescales on the order of the ion sound speed.
We also demonstrate that hybridizable discontinuous Galerkin enables solving large-scale PDE systems with implicit temporal solvers using a systematic formulation to allow automatic generation of the HDG Schur complement.
This systematic process is shown to handle over 40 coupled differential equations and to accurately compute the thousands of analytical non-zero terms of Jacobian, as well as to identify and eliminate variables which analytically vanish.
Potential future extensions include implementing bidirectional upwinding numerical fluxes\cite{cockburn2009} to improve convergence of parabolic/elliptical terms, as well as investigate effective techniques for reducing the time it takes to construct and solve the HDG Schur complement system.

\section*{Supplementary Material}

A copy of the code can be obtained from \url{https://bitbucket.org/helloworld922/dg}.

\bibliographystyle{siamplain}
\bibliography{refs.bib}
\end{document}